\documentclass[12pt,a4paper]{amsart}
\usepackage{graphicx}
\usepackage{mathrsfs}

\usepackage{enumerate}
\let\cal\mathcal

\newtheorem{theorem}{Theorem}[section]
\newtheorem{proposition}[theorem]{Proposition}
\newtheorem{lemma}[theorem]{Lemma}
\newtheorem{corollary}[theorem]{Corollary}
\newtheorem{paso}{Step}
\theoremstyle{definition}
\newtheorem{remark}[theorem]{Remark}

\newtheorem{definition}[theorem]{Definition}

\newtheorem{notation}[theorem]{Notation}

\newtheorem{example}[theorem]{Example}

\newtheorem{notation-remark}[theorem]{Notation and Remark}

\newcommand\rightmap[1]{\smash{\mathop{\rightarrow}\limits^{#1}}}

\newcommand\downmap[1]{\downarrow\rlap{$\vcenter{\hbox{$\scriptstyle#1$}}$}}
\newcommand\swmap[1]{\swarrow\rlap{$\vcenter{\hbox{$\scriptstyle#1$}}$}}

\newcommand\inj{\hookrightarrow}
\newcommand\injmap[1]{\smash{\mathop{\hookrightarrow}\limits^{#1}}}
\newcommand\surj{\mathbin{\rightarrow\mkern-6mu\rightarrow}}
\newcommand\surjmap[1]{\smash{\mathop{\surj}\limits^{#1}}}

\newcommand\cA{{\cal A}}

\newcommand\cC{{\cal C}}

\newcommand\CC{{\mathbb C}\,}

\newcommand\ZZ{{\mathbb Z}}
\newcommand\NN{{\mathbb N}}

\newcommand\PP{{\mathbb P}}
\newcommand\TT{{\mathbb T}}

\newcommand\KK{{\mathbb K}}

\newcommand\DD{{\mathbb D}}
\newcommand\AAA{{\mathbb A}}
\newcommand\fm{{\mathfrak m}}

\newcommand\mult{{\rm mult}}

\newcommand\gs{{\sigma}}
\newcommand\Char{{\rm Char}}

\newcommand\Sing{{\rm Sing}}
\newcommand\Reg{{\rm Reg}}
\newcommand\Supp{{\rm Supp}}

\newcommand\Spec{{\rm Spec}}

\newcommand\Null{{\rm Null}} 

 \newcommand\ba{{\mathbb A}}
\newcommand\bp{{\mathbb P}} 
 
 \newcommand\bff{{\mathbb F}}

\newcommand\bdm{\begin{displaymath}}
\newcommand\edm{\end{displaymath}}

\title[Essential coordinate components]
{Essential coordinate components of characteristic varieties}
\author[E. Artal]{Enrique ARTAL BARTOLO$^*$}
\address{
Departamento de Matem\'aticas\\
Universidad de Zaragoza\\
Campus Plaza San Francisco s/n\\
E-50009 Zaragoza SPAIN}
\email{artal@posta.unizar.es}
\thanks{${}^*$Partially supported by
BFM2001-1488-C02-02}
\author[J. Carmona]{Jorge CARMONA RUBER${}^\dag$}
\address{Departamento de Sistemas inform\'aticos y programaci\'on\\
Universidad Complu\-tense\\ 
Ciudad Universitaria s/n\\
E-28040 Madrid SPAIN}
\email{jcarmona@eucmos.sim.ucm.es}
\thanks{${}^\dag$ Partially supported by
BFM2001-1488-C02-01}

\author[J.I. Cogolludo]{Jos\'e Ignacio COGOLLUDO AGUST\'IN$^*$}
\address{Departamento de Matem\'aticas\\
Universidad de Zaragoza\\
Campus Plaza San Francisco s/n\\
E-50009 Zaragoza SPAIN}
\email{jicogo@posta.unizar.es}
\date\today

\begin{document}
\begin{abstract}
In this note we give an algebraic and topological
interpretation of essential coordinate components of 
characteristic varieties and illustrate their importance 
with an example. 
\end{abstract}
\maketitle

\bibliographystyle{amsplain}

\section{Introduction}
\label{section-introduction}
Let $\cC=\cC_0 \cup \cC_1 \cup \dots \cup \cC_r$ be a projective
algebraic curve in $\PP^2$. For convenience, we assume
$\cC_0$ to be a line intersecting the curve
$\cC':=\cC_1 \cup \dots \cup \cC_r$ transversally. In this
context, $\CC^2$ will refer to $\PP^2 \setminus \cC_0$. Let 
$X:=\PP^2\setminus\cC$ denote its complement and $\tilde X$
the universal abelian cover of $X$. The group of covering
transformations of  $\tilde X \rightarrow X$ acts on
$H_1(\tilde X;\ZZ)$ and endows it with a 
$\Lambda_\cC^\ZZ$-module structure, where
$\Lambda_\cC^\ZZ:=\ZZ[H_1(X;\ZZ)]$.
This module is called the Alexander module of $\cC$ and will
be denoted by $M_\cC^\ZZ$. Tensoring by a field $\KK$ one
can regard $M_\cC^\KK:=H_1(\tilde X;\KK)$ as a 
$\Lambda_\cC^\KK$-module, where $\Lambda_\cC^\KK:=\KK[H_1(X;\ZZ)]$.
When refering to $\KK=\CC$ we may drop the superscripts.

The sequence of characteristic varieties is defined in a
very natural way as a sequence of invariants of $M_\cC$ by 
means of the reduced support of its successive exterior
powers, that is,
$$\Char_k(\cC):=(\Supp_{\Lambda_\cC}(\wedge^k M_\cC))_{\text{red}}.$$
Hence each $\Char_k(\cC)$ might be regarded as a variety
in the complex torus $\Spec \Lambda_\cC = (\CC^*)^r$
whose coordinates $t_1,\dots,t_r$ are given by the 
meridians of the curve~$\cC'$ --\,see 
\cite{Libgober-characteristic} or \cite{ji-tesis} for a more 
detailed definition.

Characteristic varieties may also be considered over an
arbitrary field $\KK$, --\,cf. \cite{Matei-thesis},
\cite{Suciu-fundamental}. Such varieties, defined over the
torus $(\KK^*)^r=\Spec \Lambda^\KK_\cC$, will be denoted 
by~$\Char_k(\cC,\KK)$.

In \cite{Arapura-geometry}, Arapura proved, in a more general
context, that $\Char_k(\cC)$ consists of a finite union of tori 
translated by torsion points. 

We will define the terms {\em coordinate} and {\em essential}
when refering to an irreducible component of $\Char_k(\cC)$.

\begin{definition}
An irreducible component $V$ of $\Char_k(\cC)$ is called a
{\em coordinate component} if $V$ is contained in a 
coordinate torus, that is, if 
$V \subset \TT_i:=\{(t_1,...,t_r) \in (\CC^*)^r \mid t_i=1\}$
for some $i$.
\end{definition}

Note that the inclusion of topological spaces 
$X \inj X_i$ ($i\in \{1,...,r\}$), 
where $X_i:=\PP^2\setminus \cC(i)$ and
$\cC(i):=\cC_0 \cup \cC_1 \cup ... \cup \cC_{i-1} \cup \cC_{i+1} \cup ... \cup \cC_r$,
produces an injection of characteristic varieties as follows:
\begin{equation}
\label{eq-injection}
\Char_k(\cC(i)) \inj \Char_k(\cC).
\end{equation}

\begin{definition}
An irreducible component $V$ of $\Char_k(\cC)$ is called a
{\em non-essential component} if $V$ is contained in the image 
by (\ref{eq-injection}) for some $i$, that is, if 
$V \subset \Char_k(\cC(i))$.
Otherwise $V$ is called {\em essential}.
\end{definition}

\begin{remark}
{\rm The concepts of essential and coordinate components
depend on the embedding of $\cC$ in $\PP^2$ and not just on
its complement $X$. For instance, consider the affine
Ceva(2) arrangement 
--\,cf. 
\cite{Barthel-Hirzebruch-Hoefer-geradenkonfigurationen}
p.~80\,--
defined by $\cC_1:=(x^2-y^2)(x^2-z^2)(y^2-z^2)z$
and a conic, $Q=3(xy+xz+yz)+(x^2+y^2+z^2)$, passing through
the double points of $\cC_1$. Let $\cC=\cC_1 \cup Q$.  
Figure~\ref{ceva-q} shows the combinatorics of the 
curve~$\cC$.

\begin{figure}[ht]
\begin{center}
\hspace*{-2cm}
\includegraphics[scale=.8]{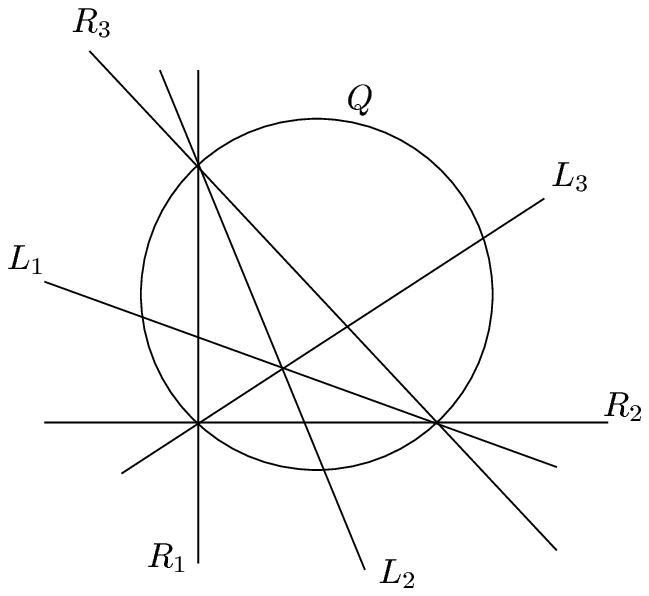}
\caption{Curve $\cC$.}
\end{center}
\label{ceva-q}
\end{figure}

Consider $r_1$, $r_2$, $r_3$,
$\ell_1$, $\ell_2$, $\ell_3$, $\ell_{\infty}$ 
and $q$ meridians around the components 
$R_1,R_2,R_3,L_1,L_2,L_3,L_{\infty}$ and $Q$ respectively. 
Hence
$$H_1(\PP^2 \setminus \cC;\ZZ)=
\langle r_i,\ell_j,q : i,j \in \{1,2,3\} \rangle 
\cong \ZZ^7.$$
Performing a Cremona transformation based on the singular
points of order four, one obtains an analytic isomorphism
of the complements which induces the following isomorphism on
homology.
$$\array{cccr}
H_1(\PP^2 \setminus \cC;\ZZ) = 
\langle r_i,\ell_j,q \rangle 
& \rightmap{\varphi} & H_1(\PP^2 \setminus \cC;\ZZ)=
\langle e_i,\ell_j,\ell_{\infty} \rangle \\
r_i & \mapsto & \ell_{\infty} - (\ell_i + e_j + e_k) \\
\ell_i & \mapsto & \ell_i \\
q & \mapsto & (e_1 + e_2 + e_3) - (\ell_1 + \ell_2 + \ell_3), \\
\endarray$$
where $e_i$ is a meridian around the exceptional component
$E_i$ resulting from blowing up the point $R_j \cap R_k$ 
($\{i,j,k\}=\{1,2,3\}$).

One can calculate $\Char_1(\cC)$ as
$${\mathbb T} = \{t_q=1\} \cap \Big(
\{t_{\ell_1}t_{\ell_2}t_{\ell_3}=t_{r_1}=t_{r_2}=t_{r_3}=1\} 
\ \ \cup$$
$$\bigcup_{\{i,j,k\}=\{1,2,3\}} 
\{t_{r_i}t_{r_j}t_{\ell_k}=t_{\ell_i}=t_{\ell_j}=t_{r_k}=1\}
\Big)$$
resulting in four (non-essential) coordinate components.

Meanwhile,
$$\varphi({\mathbb T})=
\{t_{e_1}t_{\ell_2}t_{e_3}=t_{\ell_1}t_{e_2}t_{\ell_3}\}
\ \ \cap $$
$$\Big(
\{t_{\ell_1}t_{\ell_2}t_{\ell_3}=1,t_{\infty}=
t_{\ell_i}t_{e_j}t_{e_k}\}_{\{i,j,k\}=\{i,j,k\}}
\ \ \cup$$
$$\bigcup_{\{i,j,k\}=\{1,2,3\}}
\{t_{e_i}^2=t_{\infty}=t_{\ell_i}t_{e_j}t_{e_k},
t_{\ell_j}=t_{\ell_k}=1\} \Big)$$
resulting in one (essential) non-coordinate component and
three (non-essential) coordinate components.}
\end{remark}

In the future we will refer to coordinate and
non-coordinate components of $\Char_k(\cC)$ relative to the 
torus whose coordinates are meridians of the curve
$\cC'$ in~$\CC^2$.

Essential non-coordinate components have been studied
by Libgober in~\cite{Libgober-characteristic}. In this work,
the theory of ideals of quasiadjunction was developed to 
fully determine such components. Also, an interpretation of 
the tangent cone of characteristic varieties at the identity 
was given by Cohen-Suciu~\cite{Cohen-Suciu-characteristic} 
(for hyperplane arrangements) 
and Libgober~\cite{Libgober-first} 
(for complements of curves).

Coordinate components of positive dimension are
characterized in~\cite{Libgober-characteristic} as follows.

\begin{lemma}
\label{lem-libgober}
Positive dimensional coordinate components of $\Char_1(\cC)$
are non-essential.
\end{lemma}

Our purpose in this note is to point out the importance of
essential coordinate components. In 
section~\ref{sect-coordinate} we give 
a characterization for the existence of essential 
coordinate components. In 
section~\ref{sect-zarpair} we present a Zariski pair of curves
--\,that is, with the same combinatorics but non-homeomorphic 
complements, cf.~\cite{Artal-couples}\,-- that can be 
distinguished only by essential coordinate components of their 
first characteristic varieties --\,examples of essential
coordinate components in the second characteristic varieties
were given by Cohen and Suciu 
in~\cite{Cohen-Suciu-characteristic}. In 
section~\ref{sect-ideals} we correlate essential coordinate
components to ideals of quasiadjuntion. Finally, in 
section~\ref{sect-topological} we give a topological 
interpretation of such components.

\section{Coordinate and non-essential components}
\label{sect-coordinate}
For simplicity we will use the following notation:
\begin{enumerate}

\smallbreak\item
$X_i:=\PP^2 \setminus \cC(i)$,
where 
$\cC(i):=\cC_0 \cup \dots \cup \cC_{i-1} \cup \cC_{i+1} 
\cup \dots \cup \cC_r$.

\smallbreak\item
$\Lambda^{\bullet}:=\Lambda^{\bullet}_{\cC}$,
$M^{\bullet}:=M^{\bullet}_{\cC}$,
$\Lambda^{\bullet}_i:=\Lambda^{\bullet}_{\cC(i)}$ and
$M^{\bullet}_i:=M^{\bullet}_{\cC(i)}$,
where $\bullet=\ZZ$ or $\KK$.

\smallbreak\item Note that given a ring $A$ and an $A$-module
$M$, we may define 
$$
\Char_{A,k}(M):=(\Supp_{A}(\wedge^k M))_{\text{red}}.
$$
If $A \rightmap{\varphi} B$ is a morphism of 
rings and $M$ is a $B$-module, then there is a morphism 
$\Char_{A,k}(M) \rightmap{\varphi^*} \Char_{B,k}(M)$.
In our case we will regard 
$\Char_{\Lambda_i^\KK,k}(M_i^\KK)$ 
as a subvariety of 
$\Spec \Lambda^\KK$ via the inclusion 
$\Spec {\Lambda^\KK_i} \inj \Spec \Lambda^\KK=(\KK^*)^r$,
which identifies $\Spec {\Lambda^\KK_i}$ and 
$\TT^\KK_i:=\{(t_1,...,t_r) \in (\KK^*)^r \mid t_i=1\}$. 
Such a variety, contained in $\TT^\KK_i$,
will always be denoted by $\Char_k(\cC_i,\KK)$.

\end{enumerate}

An alternative way to calculate the characteristic
varieties of a module is via the Fitting ideals of
a free resolution 
--\,cf.~\cite{Hironaka-alexander},~\cite{Matei-hall}. Let 
$$(\Lambda^\KK)^m  \ \rightmap{\phi} \ (\Lambda^\KK)^n \ 
\rightarrow M^\KK \rightarrow 0$$
be a free resolution of $M^\KK$. The map $\phi$ defines
a $n \times m$ matrix with coefficients in $\Lambda^\KK$.
The {\em Fitting ideal $F_k(M^\KK)$ of $M^\KK$} is defined as 
the ideal generated by:
$$\left\{ 
\array{ll} 
0 & {\rm if \ } k \leq \max \{0,n-m\} \\
1 & {\rm if \ } k > n \\
{\rm minors \ of \ } \phi {\rm \ of \ order \ }
(n-k+1) & {\rm otherwise.}
\endarray
\right.$$

One can therefore define
$$V_k(\cC,\KK)=V_k(M^\KK):=\Supp_{\Lambda^\KK}(M^\KK/F_k(M^\KK)).$$
An interesting remark is that 
$(V_k(\cC,\KK))_{\text{red}}=\Char_k(\cC,\KK)$ 
(cf.~\cite{Libgober-characteristic},\cite{Matei-hall}).

Suppose one has a finite presentation of $\pi_1(X)$
whose generators are meridians of the components of
$\cC' \subset \CC^2$. This allows one to obtain a finite
free resolution of the relative homology group
$M':=H_1(\tilde X,\tilde p;\ZZ)$ via Fox Calculus, 
where $\tilde p$ denotes the inverse image of
a point $p \in X$ by the universal abelian 
cover. Such a finite free resolution has the form 
\begin{equation}
\label{free-res}
(\Lambda^\KK)^m  \ \rightmap{\phi}  \ (\Lambda^\KK)^n \ 
\rightarrow M' \otimes \KK \ \rightarrow 0,
\end{equation}
where $m$ is the number of relations and $n$ the 
number of generators of the given presentation of
$\pi_1(X)$. The widely used connection between 
$\Char_k(\cC,\KK)$ and $\Char_{\Lambda^\KK,k}(M'\otimes \KK)$ 
reads as follows:
$$\Char^*_k(\cC,\KK)=\Char^*_{\Lambda^\KK,k}(M'\otimes \KK),$$
where $\Char_k^*$ denotes $\Char_k$ minus
the origin $\overline 1=(1,...,1) \in (\KK^*)^r$, 
see for instance~\cite{ji-tesis}.

\smallbreak
Consider the morphism $\sigma_i$ induced by the inclusion:
$$G:=\pi_1(X) \ \surjmap{\gs_i} \ G_i:=\pi_1(X_i).$$
This morphism is a surjection whose kernel is normally generated
by $\gamma_i$, a meridian of 
$\cC_i$. 
from the presentation of $G$ one can construct a free resolution 
of $M'_i$
$$(\Lambda^\ZZ_i)^{m+1} \  \rightmap{\phi_i} \ 
(\Lambda^\ZZ_i)^n \ \rightarrow M'_i \ \rightarrow 0$$
by adding the
relation $\gamma_i=1$. Hence, the matrix $\phi_i$ is obtained 
from $\phi$ by evaluating on $t_i=1$ and adding a column with
zeroes everywhere except in the position corresponding to 
$\gamma_i$, where a 1 is located:
\begin{equation}
\label{eq-matrix}
\phi_i=\left( \array{cccc}
\ &&\ &0\\
\ &&\ &\vdots\\
\ & &\ &0\\
\ &\phi&\ &1\\
\ & &\ &0\\
\ &&\ &\vdots\\
\ &&\ &0
\endarray
\right).
\end{equation}
Therefore one has the following:

\begin{lemma}
\label{lemma-v}
Under the previous conditions
$$\Char^*_k(\cC_i,\KK) \subset 
\Char^*_k(\cC,\KK) \cap \TT^\KK_i.$$
\end{lemma}

\begin{proof}
Using~(\ref{eq-matrix}) it is easy to check the following
inclusion of ideals 
$$F_{k+1}(M'_i) + (t_i-1) \supset 
F_{k+1}(M') + F_{k+2}(M') + (t_i-1).$$
Hence, since $F_{k+2}(M') \subset F_{k+1}(M')$, the
result follows.
\end{proof}

The inequality can be strict. This is illustrated by
Cohen-Suciu in~\cite[Example~4.4]{Cohen-Suciu-characteristic}, 
for $k=2$, using an arrangement of lines in
$\PP^2$. In this work we present an example of strict inequality
for~$k=1$. Note that, by Lemma~\ref{lem-libgober}, the
difference consists only of a finite union of points. 

\smallbreak
The variety $\Char_k(\cC,\KK) \cap \TT^\KK_i$
has the following interpretation.

\begin{lemma}
\label{lemma-intersect}
$$\Char_k(\cC,\KK) \cap \TT^\KK_i = \Char_{\Lambda^\KK,k}(M(i)^\KK),$$
where $M(i)^\KK:=M \otimes \left( \Lambda^\KK/(t_i-1) \right)$.
\end{lemma}

\begin{proof}
Tensoring a free resolution of $M^\ZZ$ by $\Lambda^\KK/(t_i-1)$ 
produces a free resolution of $M(i)^\KK$ as a
$\Lambda^\KK_i$-module. Hence, $F_k(M(i)^\KK)$ can be identified
with the ideal sum~$F_k(M^\KK) + (t_i-1)$. 
\end{proof}

Note that $M(i)^\ZZ$ is related to a certain cover of
$X$, say $\tilde X(i)$, that fits in the following
diagram
$$\array{ccc}
\tilde X & \rightmap{t_i} & \tilde X(i) \\
\downmap{\rho} & \swmap{\rho_i} & \\
X&&
\endarray$$
where $\rho$ is the universal abelian cover of $X$,
$X(i)$ is the maximal abelian cover of $X$
not ramified on $\cC_i$, and $t_i$ is the infinite 
cyclic cover of $\tilde X(i)$ ramified 
on~$\rho_i^{-1}(\cC(i))$.

The $\Lambda_i^\ZZ$-modules $M(i)^\ZZ$ and 
$H_1(\tilde X(i);\ZZ)$ fit in the following short
exact sequence

\begin{equation}
0 \to M(i)^\ZZ \to H_1(\tilde X(i);\ZZ) 
\to \Lambda_i^\ZZ/I_i=\ZZ \to 0,
\label{eqMi}
\end{equation}

where $I_i$ is the augmentation ideal of $\Lambda_i^\ZZ$.

\begin{notation}
If $\varphi$ denotes an unbranched cover of a path-connected
topological space $Z$, say $Y\rightmap{\varphi} Z$, the same 
notation will be kept for the inclusion of subgroups
$\pi_1(Y) \ \injmap{\varphi} \ \pi_1(Z)$.
\end{notation}

As a subgroup of $G$, the cover $\tilde X(i)$
corresponds to the smallest normal subgroup
$K \triangleleft G$ of $G$ containing 
$\gamma_i$ (a meridian of $\cC_i$)
such that $G/K$ is abelian. Therefore 
$K=\gs_i^{-1}(G'_i)$, where $G'_i$ denotes the 
commutator subgroup of $G_i$. Hence one has the 
following exact sequence of $\Lambda_i^\ZZ$-modules
\begin{equation}
\label{eq-ses}
0 \rightarrow \ \tilde R(i)^\ZZ \ \rightarrow \  
H_1(\tilde X(i);\ZZ) \ \rightarrow \ 
H_1(\tilde X_i;\ZZ) \ \rightarrow 0,
\end{equation}
where 
$$\tilde R(i)^\ZZ=\frac{\gs_i^{-1}(G''_i)}
{[\gs_i^{-1}(G'_i),\gs_i^{-1}(G'_i)]}.$$





Note that the morphism 
$M(i)^\ZZ \surj H_1(\tilde X_i;\ZZ)$ is
surjective.

\begin{definition}
Let us denote by $R(i)^\bullet$ the kernel
of the natural epimorphism 

$$0 \rightarrow \ R(i)^\bullet \ \rightarrow \
M(i)^\bullet \ \rightarrow \
H_1(\tilde X_i;\bullet) \ \rightarrow 0$$

obtained from~(\ref{eqMi}).

Such a $\Lambda_i^\ZZ$-module will be
called the {\em residual module of $X$ with respect to
$\cC_i$ (with coefficients in~$\bullet=\ZZ,\KK$)}.
\end{definition}

From (\ref{eqMi}) and (\ref{eq-ses}) it is obvious that
$R(i)^\ZZ$ can also be obtained as the intersection of 
$\tilde R(i)^\ZZ$ and $M(i)^\ZZ$, that is, as the kernel of

\begin{equation}
\tilde R(i)^\ZZ \to \Lambda_i^\ZZ/I_i.
\label{eq-residual}
\end{equation}

The residual module measures the existence of non-essential
coordinate components.

\begin{corollary}
The reduced residual module $R^\KK(i)_{\rm red}$ is trivial 
if and only if $\Char_1(\cC,\KK) \cap \TT^\KK_i$ consists of
non-essential components.
\end{corollary}

\begin{proof}
The result is a consequence of Lemma~\ref{lemma-intersect}
and the following remarks:
\begin{enumerate}
\item 
$M(i)^\KK=H_1(\tilde X;\ZZ) \otimes \Lambda^\KK/(t_i-1)$
\item
$M_i^\ZZ=H_1(\tilde X_i;\ZZ)$
\end{enumerate}
\end{proof}

\begin{remark} Recall that non-coordinate essential components 
may be computed by means of Libgober's method 
\cite{Libgober-characteristic}.
Then, the knowledge of non-coordinate essential components and
$R^\KK(i)_{\rm red}$, $i=1,\dots,r$ characterize essential  
coordinate components of $\Char_1(\cC,\KK)$.
\end{remark}

\section{A Zariski pair}
\label{sect-zarpair}
Consider the space ${\mathcal M}$ of sextics with the
following combinatorics: 
\begin{enumerate}
\smallbreak\item $\cC$ is a union of a smooth conic $\cC_2$
and a quartic~$\cC_4$. 
\smallbreak\item $\Sing(\cC_4)=\{P,Q\}$ where $Q$ is a cusp of
type $\AAA_4$ and $P$ is a node of type~$\AAA_1$.
\smallbreak\item $\cC_2 \cap \cC_4=\{Q,R\}$ where $Q$ is a 
$\DD_7$ on $\cC$ and $R$ is a $\AAA_{11}$ on~$\cC$.
\end{enumerate}

Performing a degenerated Cremona transformation based 
on $2 Q$ and $R$, the problem is equivalent to finding a 
nodal cubic $\tilde \cC_4$ and a smooth conic 
$\tilde \cC_2$ intersecting in two singular points of 
types $\AAA_9$ and~$\AAA_1$.

\begin{figure}[ht]
\begin{center}
\hspace*{-2cm}
\includegraphics[scale=.8]{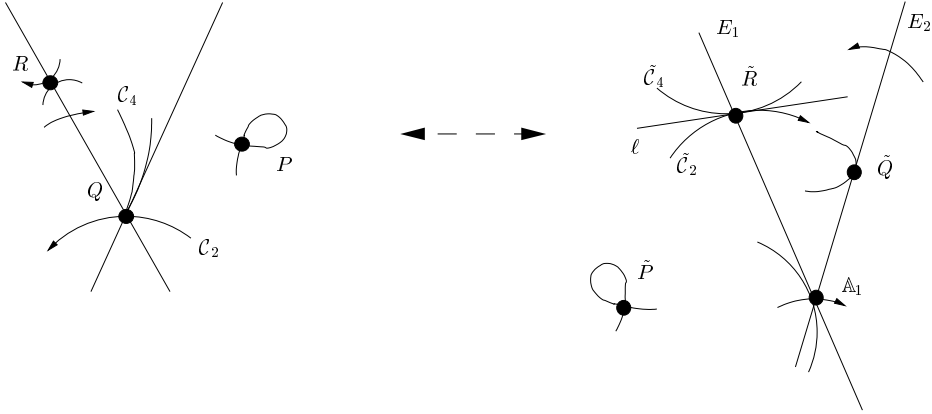}
\caption{Cremona transformation.}
\end{center}
\label{cremona}
\end{figure}

Assuming $\tilde \cC_4$ has equation $xyz+x^3-y^3$,
one can consider the following para\-metrization 
$$\array{cccc}
\varphi: & \CC & \rightarrow & \tilde \cC_4 \\
& t & \mapsto & [t:t^2:t^3-1].
\endarray$$
Note that
$\varphi|_{\CC^*}: \CC^* \rightarrow \Reg(\tilde \cC_4)$ 
is a group isomorphism from the multiplicative group
$\CC^*$ to the set of regular points on the cubic 
$\tilde \cC_4$ whose
geometric group structure has the inflexion point 
$\varphi(1)=[1:1:0]$ as unity.

Let $t_1$, $t_2$ and $t_3$ denote the parameters
corresponding to $\tilde R$, $\AAA_1$ and $\tilde Q$
respectively. One has the following relations
given by $E_1$, $\tilde \cC_2$ and~$E_2$:
$$\array{c}
t_1t_2^2=1 \\
t_1^5t_2=1 \\
t_2t_3^2=1. \\
\endarray$$
This implies that $t_1$ is a ninth  root of 1, say
$\alpha$, and $t_2=\alpha^4$. Therefore $\alpha$ must also 
be a primitive ninth root of unity. That leaves us with two 
possibilities for $t_3$, namely, $t_3=\pm\alpha^7$. The solution 
$-\alpha^7$ (resp. $+\alpha^7$) corresponds to the case where 
the tangent line to $\tilde \cC_4$ at $\tilde R$ passes 
(resp. doesn't pass) through~$\tilde Q$.

Using the techniques shown in~\cite{Artal-Carmona-ji-scwbmn},
one can obtain equations for two sextics
$\cC_6^{(1)}=\cC_4^{(1)} \cup \cC_2^{(1)}$ and
$\cC_6^{(2)}=\cC_4^{(2)} \cup \cC_2^{(2)}$ satisfying the
properties stated above and an extra property: there exists
a conic $\tilde \ell$ --\,the inverse image of $\ell$\,-- 
passing through $R$ and $Q$ such that 
$\mult_R(\tilde \ell,\cC_2^{(i)})=
\mult_R(\tilde \ell,\cC_4^{(i)})=3$,
$\mult_Q(\tilde \ell,\cC_2^{(i)})=1$ and
$\mult_Q(\tilde \ell,\cC_4^{(i)})=3+i$.

Note that, by construction, $\cC_6^{(1)}$ and
$\cC_6^{(2)}$ belong to different components of ${\mathcal M}$.
Moreover, if we consider the action of $PGL(3,\CC)$ on 
${\mathcal M}$, then ${\mathcal M}/PGL(3,\CC)$ consists of exactly
two points having representatives $\cC_6^{(1)}$ 
and~$\cC_6^{(2)}$.
\smallbreak

Special affine equations for these curves are shown below. 
The affine coordinates are $(y,z)$; the line at infinity is 
tangent to the type $\DD_7$ point, which is the base point 
of the pencil of vertical lines $y=\text{constant}$:

\begin{equation*}
f_1(y,z):= \left(  \left( y+3 \right) z+\frac{3y^2}{2} \right)
\left( {z}^{2}- \left( {y}^{2}+\frac{15}{2}\,y+\frac{9}{2} \right) 
z-3\,{y}^{3}-\frac{9y^2}{4}+\frac{y^4}{4} \right) 
\end{equation*}
for $\cC_6^{(1)}$ and 
\begin{equation*}
f_2(y,z):= \left(  \left( y+\frac{1}{3} \right) z-\frac{y^2}{6}
 \right)
\left( {z}^{2}- \left( {y}^{2}+\frac{9y}{2}+\frac{3}{2}\right) z
+\frac{y^4}{4}+\frac{3y^2}{4} \right) 
\end{equation*}
for $\cC_6^{(2)}$. 

In the future we will refer to $\cC^{(i)}$ $i=1,2$ as
the union of the sextic curve $\cC_6^{(i)}$ and a
transversal line $\cC_0$, where
$\CC^2=\PP^2 \setminus \cC_0$ 
and~$X^{(i)}=\PP^2 \setminus \cC^{(i)}$.

\begin{proposition}
\label{prop-fg}
The fundamental groups $G^{(i)}:=\pi_1(X^{(i)})$ have the 
following presentations
$$G^{(1)}=\langle e_1,e_2 : [e_2,e_1^2]=1,
(e_1e_2)^2=(e_2e_1)^2, [e_1,e_2^2]=1 \rangle$$
$$G^{(2)}=\langle e_1,e_2 : [e_2,e_1^2]=1,
(e_1e_2)^2=(e_2e_1)^2 \rangle.$$
\end{proposition}

\begin{proof}[Sketch of the proof]
The main ideas of this proof have already appeared in 
\cite{Artal-Carmona-ji-Tokunaga-sextics}.
Let us begin by considering the affine curves $\hat\cC_6^{(i)}$
defined by $f_i(y,z)=0$
and the projection $(y,z)\mapsto y$.
Note that $\CC^2\setminus\tilde\cC_6^{(i)}$
cannot be $X^{(i)}$ because of the non-generic 
choice of the line at infinity.

\begin{paso} Computation of the braid monodromy of $\hat\cC_6^{(i)}$.
\end{paso}

We follow the method introduced in 
\cite{Artal-Carmona-ji-Tokunaga-sextics} for 
``curves with real pictures'' to obtain a braid 
monodromy. In both cases, all the roots of the discriminant $\Delta$
of $f_i(y,z)$ with respect to $z$ are real.

\smallbreak
Figure~\ref{curvas} shows the real pictures of  $\hat\cC_6^{(i)}$. The
dashed lines correspond to the non-transversal vertical lines; 
the dotted stretches represent the real parts of complex conjugate solutions
and the thickened points represent real nodes with imaginary tangent lines. 
The thickened curves represent the conics.

\begin{figure}[ht]
\begin{center}
\hspace*{-2cm}
\includegraphics{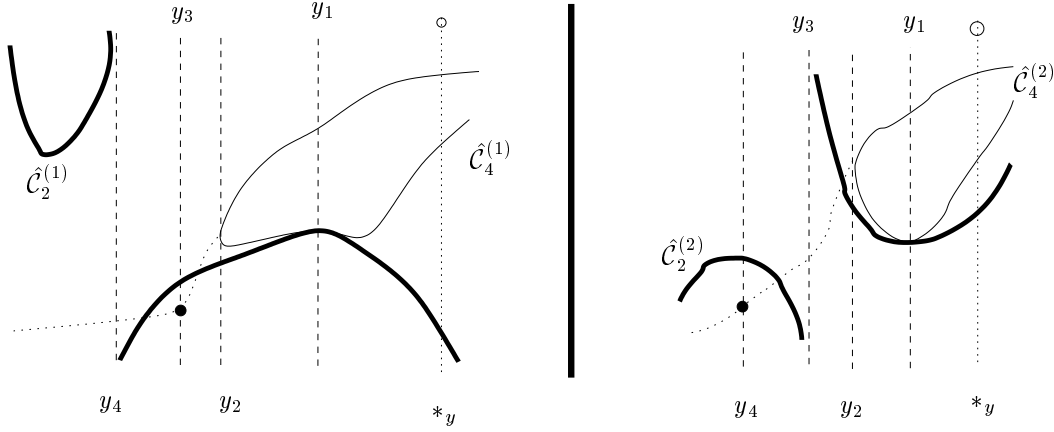}
\caption{The left figure is $\hat\cC_6^{(1)}$ and the 
right figure is $\hat\cC_6^{(2)}$.}
\end{center}
\label{curvas}
\end{figure}
 
Braid monodromy is computed for closed paths as follows.
Let us choose a \emph{big} real number $*_y$ as the base point.
We order the points in the discriminant $y_1>\dots>y_4$ and
fix $\varepsilon>0$ small enough. We use the following notation:
\begin{itemize}
\item 
$\gamma_1$ for the segment in the real line from $*_y$ to $y_1+\varepsilon$.

\item $\gamma_j$, $j=2,3,4$, for the segment in the real line from 
$y_{j-1}-\varepsilon$ to $y_{j+1}+\varepsilon$.

\item $\delta_j$, $j=1,\dots,r$ for the closed path based at 
$y_j+\varepsilon$, which runs counterclockwise along
the circle of radius $\varepsilon$ centered at $y_j$; the upper 
half is denoted by $\delta_j^+$ and the lower
half is denoted by $\delta_j^-$.
\end{itemize} 

Braid monodromy is computed along the paths 
$$
\eta_j:=\gamma_1\cdot\left(\prod_{k=2}^j(\gamma_k\cdot\delta_j^+)\right)\cdot\delta_j\cdot
\left(\gamma_1\left(\prod_{k=2}^j(\gamma_k\cdot\delta_j^+)\right)\right)^{-1},
$$
$j=1,\dots,4$. We follow the method and conventions for 
``curves with real pictures'' introduced in
\cite{Artal-Carmona-ji-Tokunaga-sextics} to obtain the braid monodromies, 
(see Table \ref{tren1} and \ref{tren2}).

\begin{table}[ht]
\begin{center}
\renewcommand\arraystretch{1.2}
\begin{tabular}{|c|c|c|c|}\hline
$\eta_1$ & $\eta_2$ &
$\eta_3$ & $\eta_4$ \\
\hline
$\sigma_2^{12}$& $\sigma_2^6*\sigma_1$ &
$\sigma_2^4*\sigma_1^2$& $(\sigma_2^4\sigma_1\sigma_2^2)*(\sigma_2\sigma_1^2\sigma_2)^{-1}$ \\
\hline
\end{tabular}
\end{center}
\caption{Braid monodromy for $\hat\cC_6^{(1)}$ }
\label{tren1}
\end{table}
\vspace*{-1cm}
\begin{table}[ht]
\begin{center}
\renewcommand\arraystretch{1.2}
\begin{tabular}{|c|c|c|c|}\hline
$\eta_1$ & $\eta_2$ &
$\eta_3$ & $\eta_4$ \\
\hline
$\sigma_2^{12}$& $\sigma_2^6*\sigma_1$ &
$(\sigma_2^5\sigma_1)*(\sigma_1\sigma_2^2\sigma_1)^{-1}$& $\sigma_2^2*\sigma_1^2$ \\
\hline
\end{tabular}
\end{center}
\caption{Braid monodromy for $\hat\cC_6^{(2)}$}
\label{tren2}
\end{table}

The notation $\sigma*\tau$ represents $\sigma\tau\sigma^{-1}$.
\begin{paso} Computation of the fundamental group of $\hat\cC_6^{(i)}$.
\end{paso}

Note that the standard Zariski-van Kampen method cannot be applied in this instance
due to vertical asymptotes. A generalized Zariski-van Kampen method applicable
to this case can be found in \cite{Artal-Carmona-ji-Luengo-Melle}. 

Let us choose
a basis $a_1,a_2,a_3$ of the fundamental group of the fiber on $*_y$ 
as in the case of the fundamental group of the complement 
of the discriminant. These elements also generate 
$\pi_1(\CC^2\setminus\hat\cC_6^{(i)})$. The standard action of the 
braid group on the free group generated by $a_1,a_2,a_3$ is defined as
$$
a_j^{\sigma_i}=
\begin{cases}
a_j&\text{ if $j\neq i,i+1$}\\
a_{j+1}&\text{ if $j=i$}\\
a_{j+1}*a_i&\text{ if $j=i+1$.}
\end{cases}
$$
Let us denote by $\tau_i$ the image of $\eta_i$ by the braid 
monodromy. If $\eta_j$ does not correspond to an asymptote, 
the relations $a_j=a_j^{\tau_i}$, $j=1,2,3$ are satisfied.
Otherwise the relations 
$a_j^{b_i}=a_j^{\tau_i}$, $j=1,2,3$ are satisfied, where
$$
b_i=
\begin{cases}
a_3^{a_2 a_3 a_2 a_1}&\text{if $i=4$ for $\tilde\cC_6^{(1)}$}\\
a_3^{a_2 a_3 a_2 a_3 a_2 a_1}&\text{if $i=3$ for $\tilde\cC_6^{(2)}$.}
\end{cases}
$$
The generalized Zariski-van Kampen method asserts that the above relations
are sufficient to give a presentation of $\pi_1(\CC^2\setminus\hat\cC_6^{(i)})$.

\begin{paso} Computation of the fundamental group of $\cC_6^{(i)}$.
\end{paso}
In order to obtain $\pi_1(\PP^2\setminus\cC_6^{(i)})$, it is enough
to factor $\pi_1(\CC^2\setminus\hat\cC_6^{(i)})$ by a meridian of the 
line at infinity. Blowing up the projection point one can see the 
fibration used for the braid monodromy
as the restriction of the standard ruled fibration $\bff_1\to\bp^1$. 

Figure~\ref{infinito} gives the real picture of $\hat\cC_6^{(i)}$ around 
both the exceptional divisor (the horizontal line) and the strict transform 
of the line at infinity (the vertical line). 
Since the self-intersection number of the exceptional curve
equals $-1$, it is easily seen that
$(b_j a_2 a_1 a_3 a_2 a_1)^{-1}$ is a meridian of the line at infinity.

\begin{figure}[ht]
\begin{center}
\hspace*{-2cm}
\includegraphics{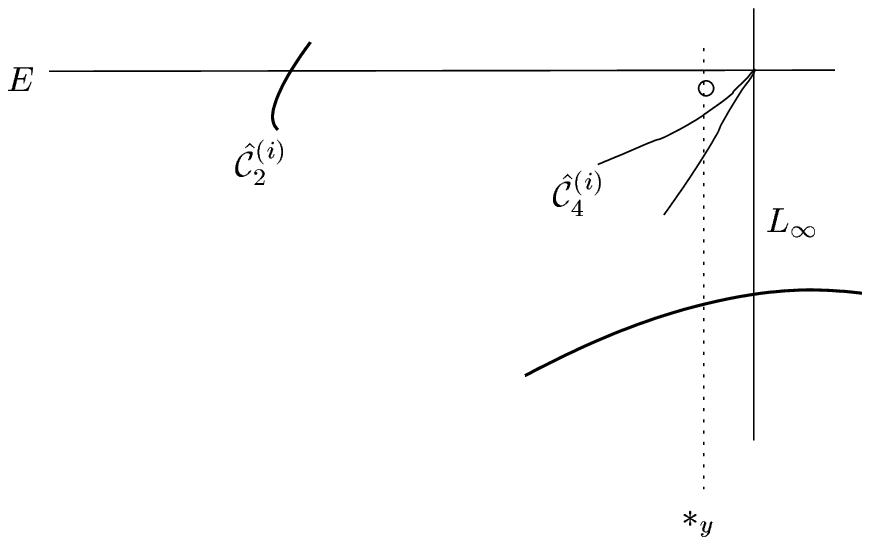}
\caption{}
\label{infinito}
\end{center}
\end{figure}

We apply GAP4 \cite{GAP4} to perform Tietze transformations
in order to obtain the following presentations of
the fundamental groups:

$$\pi_1(\cC_6^{(1)})=\langle a_1,a_3 : [a_3,a_1^2]=1,
(a_1 a_3)^2=(a_3 a_1)^2, a_3 a_1 a_3 a_1^3=1 \rangle,$$
$$\pi_1(\cC_6^{(2)})=\langle a_1,a_3 : [a_3,a_1^2]=1,
(a_1 a_3)^2=(a_3 a_1)^2, (a_3 a_1^2)^2=1 \rangle.$$

\begin{paso} Computation of the fundamental group of $\cC_6^{(i)}$.
\end{paso}

We apply last step and \cite[Prop. 2]{Artal-Carmona-ji-Tokunaga-sextics}
to obtain the presentations of Proposition~\ref{prop-fg}.
\end{proof}

One can easily calculate the matrix of the free resolution
of $M^{'(i)}$ by means of Fox Calculus to get
$$\phi^{(2)}=\left[
\array{cc}
(t_1-1)(t_1+1) & (t_1-1)(t_1t_2+1)\\
(t_2-1)(t_1+1) & (t_2-1)(t_1t_2+1)
\endarray
\right]$$

and

$$\phi^{(1)}=\left[
\array{ccc}
(t_1-1)(t_1+1) & (t_1-1)(t_1t_2+1) & (t_1-1)(t_2+1)\\
(t_2-1)(t_1+1) & (t_2-1)(t_1t_2+1) & (t_2-1)(t_2+1)
\endarray
\right].$$

Hence
$$\array{ccc}
\Char^*_1(\cC^{(2)},\KK)=
\left\{ \array{cl} \emptyset & {\rm if \ } char_\KK = 2\\
(-1,1) & {\rm otherwise} \endarray \right.
& {\rm and} &
\Char^*_1(\cC^{(1)},\KK)=\emptyset.
\endarray$$

Using Reidemeister-Schreier's method\footnote{
For computationally more difficult problems, the authors 
have developed a software package programmed for Singular 
--\,\cite{Singular}, \cite{Singular-homologlib}\,-- which
is available upon request.}
one can explicitly
obtain the short exact sequences described in~(\ref{eq-ses})
for both meridians $e_1$ and~$e_2$:

\begin{enumerate}
\item For the coordinate components on $t_1=1$:
$$
\array{ccccccccc}
0 & \rightarrow & \tilde R^{(i)\ZZ}(1) & \rightarrow &  
\ZZ e_1 \oplus \ZZ/2\ZZ [e_1,e_2]& \rightarrow & 
0 & \rightarrow & 0.
\endarray
$$

From~(\ref{eq-residual}) one has $e_1 \mapsto e_1$ and 
$[e_1,e_2] \mapsto (t_2-1)e_1$. This implies 
$R^{(1)\ZZ}(1)=R^{(2)\ZZ}(1)=\ZZ/2\ZZ[e_1,e_2]$, 
and therefore

$$R^{(1)\KK}(1)=R^{(2)\KK}(1)=\left\{ 
\array{ll} \KK & \text{if } char_\KK=2\\
0 & \text{otherwise.} \endarray \right.$$

\item For the coordinate components on $t_2=1$:
$$
\array{ccccccccc}
0 & \rightarrow & \tilde R^{(1)\ZZ}(2) & \rightarrow &
\ZZ e_2 \oplus \ZZ/2\ZZ [e_1,e_2] & \rightarrow & 0 &
\rightarrow & 0 \\
0 & \rightarrow & \tilde R^{(2)\ZZ}(2) & \rightarrow &  
\ZZ e_2 \oplus \ZZ [e_1,e_2] & \rightarrow & 0 & 
\rightarrow & 0.
\endarray
$$

From~(\ref{eq-residual}) one has $e_2 \mapsto e_2$ and 
$[e_1,e_2] \mapsto (t_1-1)e_2$. This implies
$R^{(1)\ZZ}(2)=\ZZ/2\ZZ[e_1,e_2]$ and
$R^{(2)\ZZ}(2)=\ZZ[e_1,e_2]$. Therefore

$$\array{ccc}
R^{(1)\KK}(2)=\begin{cases} \KK & \text{if } char_\KK=2\\ 
0 & {\rm otherwise} \end{cases} &
\text{and} &
R^{(2)\KK}(2)=\!\KK.
\endarray$$
\end{enumerate}

\section{Ideals of quasiadjunction}
\label{sect-ideals}
Global ideals of quasiadjunction 
--\,cf.~\cite{Libgober-characteristic}\,-- 
give the required topological information to calculate the
first Betti number of abelian covers of $\PP^2$ ramified
along a curve $\cC$. This is based on Sakuma's formula 
--\,cf. Theorem~7.3. in~\cite{Sakuma-homology}\,-- for
the first Betti number of abelian covers of compact smooth
complex surfaces ramified along complex curves. In our case 
it can be described as follows. Let the manifold 
$\overline X_{\rho}$ denote 
a certain abelian cover of $\PP^2$ ramified along $\cC$ and
let $X_{\rho}$ be its associated unbranched covering.
Consider the set $\TT_{X_{\rho}}$ of representations 
$H_1(X_\rho;\ZZ) \to \CC^*$ and 
$\TT^*_{X_\rho}:=\TT_{X_\rho} \setminus \{1\}$.
Then
\begin{equation}
\label{sakuma-branched}
b_1(\overline X_\rho)=\sum_{\xi \in \TT^*_{X_\rho}} 
\Null(\cC_{\xi};\xi),
\end{equation}
where
$$\Null(\cC_{\xi};\xi):=
\max\{k \in \NN \mid \xi \in \Char_k(\cC_{\xi})\},$$
and $\Char_k(\cC_{\xi})\subset\TT^\xi$ which is the
coordinate subtorus of $\TT_{X_\rho}$ determined by
the trivial coordinates of $\xi$. The curve $\cC_{\xi}$ 
denotes the union of components of $\cC$ associated 
with the non-trivial coordinates of $\xi$.
Libgober proved that $\Null(\cC_{\xi};\xi)$ coincides with
the irregularity of a certain ideal sheaf on $\PP^2$
associated with $\xi$ and~$\cC$.

Alternatively, Sakuma also provided the following formula
for the first Betti number of~$X_\rho$:
\begin{equation}
\label{sakuma-unbranched}
b_1(X_\rho)=b_1(X) + \sum_{\xi \in \TT^*_{X_\rho}} 
\Null(\cC;\xi).
\end{equation}

\begin{remark} Note that $\Null(\cC;\xi)=\Null(\cC_\xi;\xi)$ if
no coordinate of $\xi$ equals $1$. In general equality holds only if
$\Char_k(\cC)\cap\TT^\xi=\Char_k(\cC_{\xi})$. 
\end{remark}

Equation (\ref{sakuma-unbranched}) shows that coordinate
components have an effect on the topology of the unbranched
cover $X_\rho$.

As mentioned in the Introduction, Libgober described an
algebraic method to compute non-coordinate components. Using
deletion of components, this method also exhibits a way to
compute non-essential coordinate components. It is
still an open problem whether such method exists for
essential coordinate components or not. 
Even though the proof doesn't cover the essential coordinate
case, Libgober's method works in certain cases as the
following example shows.

\begin{example}
{\rm 
For the sake of simplicity we follow the notation 
introduced in~\cite{Libgober-characteristic} (Section~2). 
Consider the sextic curves $\cC^{(1)}$ and $\cC^{(2)}$
described in the previous section. The contributing ideals
of quasiadjunction are described as follows:
\begin{enumerate}
\item $\cA_{(0,\frac{1}{2})}(1)$. It is supported on
$\{Q,R\}$ and its stalks at these points are 
$\fm_1$ and $\fm_3$ respectively, where
$\fm_k=(x^k,y)$ and $(x,y)$ denote a local system of 
coordinates for which the equation of the singular point
$\AAA_n$ is~$y^2-x^{n+1}$.
\item $\cA_{(\frac{1}{2},0)}(2)$. It is supported on
$\{Q,R\}$ and its stalk at both points is~$\fm_3$.
\item $\cA_{(\frac{2}{3},\frac{1}{6})}(1)$. It is supported
on $\{Q,R\}$ and its stalks at these points are 
$\fm_2$ and $\fm_1$ respectively.
\end{enumerate}

The only ideal sheaf with non-trivial irregularity is
$\cA_{(\frac{1}{2},0)}(2)$ relative to $\cC^{(2)}$, due to
the existence of the conic $\tilde \ell$ passing through
three infinitely near points of $Q$ and $R$. As in the
essential case, this irregularity describes the existence of
a non-trivial point in $\Char_1(\cC^{(2)})$, namely the
point~$\exp(\frac{1}{2},0)=(-1,1)$.}
\end{example}

\section{A topological approach}
\label{sect-topological}
Finally, we will describe the somewhat surprising difference
in the first Betti numbers of unbranched covers of the
complement of sextics whose associated branched cover of
$\PP^2$ only ramifies on a conic. First let us give the 
commutator exact sequences for $G^{(1)}$ and $G^{(2)}$:
$$0 \ \rightarrow [G^{(i)},G^{(i)}]\ \rightarrow \
G^{(i)} \ \rightarrow \ZZ e_1 \oplus \ZZ e_2 \ \rightarrow \ 0.$$
They were obtained using Reidemeister-Schreier's method.
Note that $[G^{(i)},G^{(i)}]$ is a cyclic group generated
by $t:=[e_1,e_2]$ and hence
$[G^{(i)},G^{(i)}]=M^{(i)\ZZ}$. The computation gives
that  $M^{(1)\ZZ}$ is the cyclic group of
order $2$, whereas $M^{(2)\ZZ}$ is the infinite cyclic group.
Moreover, the $\Lambda^{\ZZ}$-module structure of $M^{(i)\ZZ}$ is
given as follows:
$$\array{c}
t_1[e_1,e_2]=e_1[e_1,e_2]e_1^{-1}=[e_1,e_2]^{-1}=[e_1,e_2]
\Rightarrow t_1 \cdot t=t\\
\\
t_2[e_1,e_2]=e_2[e_1,e_2]e_2^{-1}=[e_1,e_2] \Rightarrow t_2 \cdot t=t
\endarray$$
in $M^{(1)\ZZ}$ and 
$$\array{c}
t_1[e_1,e_2]=e_1[e_1,e_2]e_1^{-1}=[e_1,e_2]^{-1} \Rightarrow t_1 \cdot t=t^{-1}\\
\\
t_2[e_1,e_2]=e_2[e_1,e_2]e_2^{-1}=[e_1,e_2] \Rightarrow t_2 \cdot t=t
\endarray$$
in~$M^{(2)\ZZ}$. 
This gives the corresponding characteristic varieties over
any field. 

According to Sakuma's formula in the unbranched 
case~(\ref{sakuma-unbranched}), we should look into 
the branched covering 
$\rho^{(i)}:\overline X_{\rho}^{(i)}\to\PP^2$
ramified along the smooth conic $\cC_2^{(i)}$ in
$\cC^{(i)}$. It is well known that $\overline X_{\rho}^{(i)}$
is isomorphic to $\PP^1\times\PP^1$ and the preimage of the
conic corresponds to the diagonal curve. A simple local
computation shows that the preimage of the quartic curve is
a curve of bidegree $(4,4)$ having one $\AAA_9$, one~$\AAA_5$
and two $\AAA_1$ points. One can perform a standard Cremona 
transformation from $\PP^1\times\PP^1$ to $\PP^2$ which sends 
$(\rho^{(i)})^{-1}(\cC_4^{(i)})$ to a sextic curve having 
four nodes, two $\ba_3$ points and one $\ba_7$ point.
Such curves must have two irreducible components.

\smallbreak
An alternative way to calculate the first Betti number of 
$X_{\rho}^{(i)}$, the unbranched cover associated with
$\rho^{(i)}$, is by direct calculation of the groups. One
can compute the fundamental group of $X_{\rho}^{(i)}$ by
means of Reidemeister-Schreier's method and thus the rank of
its abelianization provides the required Betti number.

Since the Betti number of the complement of curves in
$\PP^1\times\PP^1$ can be computed from the bidegrees of their
irreducible components, one deduces that the two irreducible 
components of $(\rho^{(2)})^{-1}(\cC_4^{(2)})$ have bidegrees 
$(2,2)$ whereas those of $(\rho^{(1)})^{-1}(\cC_4^{(1)})$
have bidegrees $(3,1)$ and $(1,3)$.

\providecommand{\bysame}{\leavevmode\hbox to3em{\hrulefill}\thinspace}
\providecommand{\MR}{\relax\ifhmode\unskip\space\fi MR }
\providecommand{\MRhref}[2]{%
  \href{http://www.ams.org/mathscinet-getitem?mr=#1}{#2}
}
\providecommand{\href}[2]{#2}

\end{document}